\documentclass{article}

\usepackage{amsmath,amsthm,bm}
\usepackage{authblk}
\usepackage{enumerate}
\usepackage{comment}
\usepackage{amssymb}
\usepackage{pgfplots}
\usepackage{subcaption}
\usepackage{graphicx}
\usepackage{url}
\usetikzlibrary{arrows,positioning,calc,patterns}

\newtheorem{thm}{Theorem}[section]
\newtheorem{cor}[thm]{Corollary}
\newtheorem{conj}[thm]{Conjecture}

\newtheorem{lem}[thm]{Lemma}

\makeatletter
\newcommand*{\house}[1]{%
  \mathord{%
    \mathpalette\@house{#1}%
  }%
}
\newcommand*{\@house}[2]{%
  \dimen@=\fontdimen8 %
      \ifx#1\scriptscriptstyle\scriptscriptfont
      \else\ifx#1\scriptstyle\scriptfont
      \else\textfont\fi\fi
      3 %
  \sbox0{%
    $#1%
      \vrule width\dimen@\relax
      \overline{%
        \kern2\dimen@
        \begingroup % to keep changes of \dimen@ in #2 local
          #2%
        \endgroup
        \kern2\dimen@
      }%
      \vrule width\dimen@\relax
      \mathsurround=1.5\dimen@ % outside margin
    $%
  }%
  \ht0=\dimexpr\ht0-\dimen@\relax
  \dp0=\dimexpr\dp0+2\dimen@\relax
  \vbox{%
    \kern\dimen@ % reinsert previously removed space
    \copy0 %
  }%
}

\newcommand{\R}{\mathbb R}
\newcommand{\Z}{\mathbb Z}
\newcommand{\Q}{\mathbb Q}

\DeclareMathOperator{\Aut}{\operatorname{Aut}}
\DeclareMathOperator{\orb}{\operatorname{orb}}

\begin{document}

\title{Absolute Bound On the Number of Solutions of Certain Diophantine Equations of Thue and Thue-Mahler Type}

\author{Anton Mosunov}
\affil{University of Waterloo}

%%% MATHEMATICS SUBJECT CLASSIFICATION 11D59 %%%

\date{}

\maketitle

\begin{abstract}
Let $F \in \Z[x, y]$ be an irreducible binary form of degree $d \geq 7$ and content one. Let $\alpha$ be a root of $F(x, 1)$ and assume that the field extension $\Q(\alpha)/\Q$ is Galois. We prove that, for every sufficiently large prime power $p^k$, the number of solutions to the Diophantine equation of Thue type
$$
|F(x, y)| = tp^k
$$
in integers $(x, y, t)$ such that $\gcd(x, y) = 1$ and $1 \leq t \leq (p^k)^\lambda$ does not exceed $24$. Here $\lambda = \lambda(d)$ is a certain positive, monotonously increasing function  that approaches one as $d$ tends to infinity. We also prove that, for every sufficiently large prime number $p$, the number of solutions to the Diophantine equation of Thue-Mahler type
$$
|F(x, y)| = tp^z
$$
in integers $(x, y, z, t)$ such that $\gcd(x, y) = 1$, $z \geq 1$ and $1 \leq t \leq (p^z)^{\frac{10d - 61}{20d + 40}}$ does not exceed 1992. Our proofs follow from the combination of two  principles of Diophantine approximation, namely the generalized non-Archimedean gap principle and the Thue-Siegel principle.
\end{abstract}

\section{Introduction}

In this article we analyze certain Diophantine equations of Thue and Thue-Mahler type. A \emph{Thue equation} is an equation of the form
\begin{equation} \label{eq:thue_equation}
F(x, y) = m,
\end{equation}
where $F \in \Z[x, y]$ is a homogeneous polynomial of degree $d \geq 3$ with nonzero discriminant $D(F)$, $m$ is a fixed positive integer, and $x, y$ are integer variables. In 1909 Thue \cite{thue} established that there is a finite upper bound on the number of solutions of (\ref{eq:thue_equation}), provided that $F$ is irreducible. Since Thue's time, the estimates on the number of solutions of (\ref{eq:thue_equation}) have been improved significantly. In 1933, assuming that $F$ is irreducible, Mahler established the existence of a number $C$, dependent only on $F$, such that the number of primitive solutions to (\ref{eq:thue_equation}), --- that is, solutions with $x$ and $y$ coprime, --- does not exceed $C^{1 + \omega(m)}$, where $\omega(m)$ denotes the number of distinct prime divisors of $m$ \cite{mahler2}. In fact, his result was even stronger: if instead of (\ref{eq:thue_equation}) we consider the equation
\begin{equation} \label{eq:thue-mahler_equation}
F(x, y) = p_1^{k_1}\cdots p_t^{k_t},
\end{equation}
where $p_1, p_2, \ldots, p_t$ are distinct fixed prime numbers, then it follows from Mahler's argument that the number of integer solutions $(x, y, k_1, k_2, \ldots, k_t)$ to (\ref{eq:thue-mahler_equation}), with $x, y$ coprime and $k_i$ non-negative, does not exceed $C^{1 + t}$. The equation (\ref{eq:thue-mahler_equation}) is called a \emph{Thue-Mahler equation}. Further improvements to this estimate have been made by Erd\H{o}s and Mahler \cite{erdos-mahler}, and Lewis and Mahler \cite{lewis-mahler}.

It was conjectured by Siegel that the number of primitive solutions to (\ref{eq:thue_equation}) should not depend on the coefficients of $F$. Siegel's conjecture was established in 1984 by Evertse \cite{evertse2}, who proved that the number of primitive solutions to (\ref{eq:thue-mahler_equation}) does not exceed
\begin{equation} \label{eq:evertses_estimate}
2\cdot 7^{d^3(2t + 3)},
\end{equation}
where a binary form $F$ of degree $d$ was assumed to be divisible by at least three pairwise linearly independent linear forms in some algebraic number field. An estimate on the number of solutions to (\ref{eq:thue_equation}) thus follows by replacing the number $t$ in (\ref{eq:evertses_estimate}) with $\omega(m)$.

When integers $x$ and $y$ are arbitrary, the number of solutions to (\ref{eq:thue_equation}) can be large. For example, in 2008 Stewart \cite{stewart3} proved that when $F$ is of degree $3$ and $D(F) \neq 0$ then there is a positive number $c = c(F)$ such that the number of solutions to (\ref{eq:thue_equation}) is at least $c(\log m)^{1/2}$. However, if we restrict our attention only to primitive solutions, then their number does not seem to increase with the growth of $m$. In 1987 it was conjectured by Erd\H{o}s, Stewart and Tijdeman \cite{erdos-stewart-tijdeman} that the number of primitive solutions to (\ref{eq:thue_equation}) does not exceed some constant, which depends only on $d$. In the same year Bombieri and Schmidt \cite{bombieri-schmidt} proved that the number of primitive solutions to (\ref{eq:thue_equation}) does not exceed
$$
Cd^{1 + \omega(m)},
$$
where the constant $C$ is absolute. In 1991 Stewart \cite{stewart2} replaced $\omega(m)$ in the above estimate with $\omega(g)$, where $g$ is a divisor of $m$ satisfying $g \gg_F m^{(4+d)/3d}$ (this is the statement of \cite[Theorem 1]{stewart2} with $\varepsilon = 1/2$). In the same paper, Stewart conjectured the following.

%The following three conjectures associated with (\ref{eq:thue_equation}) remain open.
%
%\textbf{Schmidt's Conjecture.} The number of solutions to (\ref{eq:thue_equation}) does not exceed $C(\log m)^c$, where $c$ is absolute and $C = C(F)$ \cite[Chapter III]{schmidt3}.
%
%\textbf{Erd\"os-Stewart-Tijdeman Conjecture.} There exists a number $C = C(d)$ such that equation (\ref{eq:thue_equation}) has at most $C$ solutions in coprime integers $x, y$ \cite{erdos-stewart-tijdeman}.

\begin{conj} \label{conj:stewart}
\emph{(Stewart, \cite[Section 6]{stewart2})} There exists an absolute constant $c_0$ such that for any binary form $F \in \Z[x, y]$ with nonzero discriminant and degree at least three there exists a number $C = C(F)$, such that if $m$ is an integer larger than $C$, then the Thue equation (\ref{eq:thue_equation}) has at most $c_0$ solutions in coprime integers $x$ and $y$.
\end{conj}

%Note that the condition $m > C(F)$ in the above conjecture cannot be removed due to a counterexample
%%
%$$
%\left|x^d + a(x - y)(2x - y)\cdots (dx - y)\right| = 1.
%$$
%%
%As remarked by Schmidt \cite[Chapter III]{schmidt3}, this equation has at least $2d$ primitive solutions of the form $(\pm 1, j)$, where $j \in \{-d, -d + 1, \ldots, d\}$. Further, when $a$ is prime then the above polynomial is irreducible over $\Q$ by Eisenstein's criterion, so certainly it has a nonzero discriminant. Thus, at least when $m$ is small, the number of solutions can depend on $d$.\footnote{However, it was remarked by Erd\H{o}s, Stewart and Tijdeman \cite{erdos-stewart-tijdeman} that the number of primitive solutions to (\ref{eq:thue_equation}) could be bounded above by a constant, which only on $d$.}

The most notable step forward towards Conjecture \ref{conj:stewart} can be found in the work of Thunder \cite{thunder2}. Based on \cite{stewart2} he gives a heuristic that supports the conjecture of Stewart when the degree of the form $F$ is at least five. By using a generalization of the non-Archmiedean gap principle established in \cite{mosunov2}, we develop new methods for estimating the number of primitive solutions of (\ref{eq:thue_equation}) and (\ref{eq:thue-mahler_equation}) in the case $t = 1$, thus providing theoretical evidence in support of Stewart's conjecture. Instead of looking at (\ref{eq:thue_equation}) and (\ref{eq:thue-mahler_equation}) though, we study equations of the form
$$
|F(x, y)| = tp^z,
$$
with $p^z$ a prime power and $t$ an integer variable, which is ``small'' in comparison to $p^z$. We demonstrate that it is possible to provide an absolute bound on the number of primitive solutions, provided that $F$ is irreducible of degree $d \geq 7$ and the order of the Galois group of $F(x, 1)$ over $\Q$ is equal to $d$.

In order to state the main results given in Theorems \ref{thm:main_theorem_thue} and \ref{thm:main_theorem_thue-mahler}, we need to introduce the notion of an \emph{enhanced automorphism group} of a binary form. For a $2 \times 2$ matrix $M = \left(\begin{smallmatrix}s & u\\t & v\end{smallmatrix}\right)$, with complex entries, define the binary form $F_M$ by
$$
F_M(x, y) = F(sx + uy, tx + vy).
$$
Let $\overline{\Q}$ denote the algebraic closure of the rationals and let $K$ be a field containing $\mathbb Q$. We say that a matrix $M = \left(\begin{smallmatrix}s & u\\t & v\end{smallmatrix}\right) \in \operatorname{M}_2(K)$ is a \emph{$K$-automorphism of $F$} (resp., $|F|$) if $F_M = F$ (resp., $F_M = \pm F$). The set of all $K$-automorphisms of $F$ (resp., $|F|$) is denoted by $\Aut_K F$ (resp., $\operatorname{Aut}_K |F|$). We define
\begin{equation} \label{eq:G}
\Aut' |F| = \left\{\frac{1}{\sqrt{|sv - tu|}}\begin{pmatrix}s & u\\t & v\end{pmatrix} \colon s, t, u, v \in \Z\right\} \cap \Aut_{\overline{\Q}} |F|
\end{equation}
and refer to it as the \emph{enhanced automorphism group} of $F$. See \cite[Lemma 7.2]{mosunov2} for a proof that $\Aut' |F|$ contains at most $24$ elements, provided that $d \geq 3$.

For a nonzero polynomial $P \in \mathbb Z[x_1, \ldots, x_n]$, we define the \emph{content} of $P$ to be the greatest common divisor of its coefficients. For an arbitrary finite set $X$, let $\#X$ denote its cardinality. Let
\begin{equation} \label{eq:bound_on_lambda}
f(d) = \frac{20d - 41}{80}\left(\frac{\sqrt{d^2 + 16d}}{d} - 1\right) - 1
\end{equation}
and notice that $f(d)$ is a positive monotonously increasing function on the interval $[7, \infty)$, which approaches one as $d$ tends to infinity. We prove the following.

\begin{thm} \label{thm:main_theorem_thue}
Let $F \in \Z[x, y]$ be an irreducible binary form of degree $d \geq 7$ and content one. Let $\alpha$ be a root of $F(x, 1)$ and assume that the field extension $\Q(\alpha)/\Q$ is Galois. Let $\lambda$ be such that $0 \leq \lambda < f(d)$, where $f(d)$ is defined in (\ref{eq:bound_on_lambda}). Let $p$ be prime, $k$ a positive integer, and consider the Diophantine equation
\begin{equation} \label{eq:thue_equation_1}
|F(x, y)| = tp^k.
\end{equation}
Provided that $p^k$ is sufficiently large, the number of solutions to (\ref{eq:thue_equation_1})  in integers $(x, y, t)$ such that
$$
\gcd(x, y) = 1 \quad \text{and} \quad 1 \leq t \leq (p^k)^\lambda
$$
is at most $\#\Aut' |F|$. In particular, it does not exceed $24$. More precisely, for any two solutions $(x_1, y_1, t_1), (x_2, y_2, t_2)$ there exists a matrix $M = |sv - tu|^{-1/2}\cdot \left(\begin{smallmatrix}s & u\\t & v\end{smallmatrix}\right)$ in $\Aut' |F|$ such that
$$
\frac{x_2}{y_2} = \frac{sx_1 + uy_1}{tx_1 + vy_1}.
$$
\end{thm}

\begin{thm} \label{thm:main_theorem_thue-mahler}
Let $F \in \Z[x, y]$ be an irreducible binary form of degree $d \geq 7$ and content one. Let $\alpha$ be a root of $F(x, 1)$ and assume that the field extension $\Q(\alpha)/\Q$ is Galois. Let $\lambda$ be such that
$$
0 \leq \lambda < 1 - 8.1/(d + 2).
$$
Let $p$ be prime and consider the Diophantine equation
\begin{equation} \label{eq:thue-mahler_equation_1}
|F(x, y)| = tp^z.
\end{equation}
Provided that $p$ is sufficiently large, the number of solutions to (\ref{eq:thue-mahler_equation_1}) in integers $(x, y, z, t)$ such that
$$
\gcd(x, y) = 1, \quad z \geq 1 \quad \text{and} \quad 1 \leq t \leq (p^z)^\lambda
$$
is at most
$$
\#\Aut' |F|\cdot \left\lfloor 1 + \frac{11.51 + 1.5\log d + \log\left((d - 2.05)/(1 + \lambda)\right)}{\log((d - 2.05)/(1 + \lambda) - 0.5d)}\right\rfloor.
$$
\end{thm}

If we let $\lambda(d) = 0.5 - 4.05/(d + 2)$, then the function
$$
g(d) = 1 + \frac{11.51 + 1.5\log d + \log\left((d - 2.05)/(1 + \lambda(d))\right)}{\log((d - 2.05)/(1 + \lambda(d)) - 0.5d)}
$$
is monotonously decreasing on the interval $[7, \infty)$. Since $g(7) \approx 83.3$, we can use the upper bound $\#\Aut'|F| \leq 24$ as well as Theorem \ref{thm:main_theorem_thue-mahler} to conclude that the number of solutions $(x, y, z, t)$ to (\ref{eq:thue-mahler_equation_1}) satisfying the aforementioned conditions does not exceed $24 \cdot \lfloor g(7)\rfloor = 1992$ when $d \geq 7$. Furthermore, since $g(10^{15}) < 4$ and $\lim\limits_{d \rightarrow \infty} g(d) = 3.5$, we can also conclude that it does not exceed $24 \cdot \lfloor g(10^{15})\rfloor = 72$ when $d \geq 10^{15}$.

The proof of Theorem \ref{thm:main_theorem_thue} follows from the generalized non-Archimedean gap principle, whose statement is given in Section \ref{sec:theory}, Lemma \ref{lem:non-archimedean_gap_principle}. The proof of Theorem \ref{thm:main_theorem_thue-mahler} follows from the combination of the non-Archimedean gap principle and the Thue-Siegel principle, as formulated by Bombieri and Mueller \cite{bombieri-mueller}. Both of these principles have been utilized in \cite{mosunov2} so to establish the result stated in Lemma \ref{lem:count_large_approximations}. Unfortunately, due to the application of Roth's Theorem \cite{roth} it is not yet possible to determine how large a prime power $p^k$ in Theorem \ref{thm:main_theorem_thue} or a prime $p$ in Theorem \ref{thm:main_theorem_thue-mahler} should be in order for the respective absolute bound to hold. The author expects that it is possible to overcome this problem if one is able to generalize the non-Archimedean gap principle even further and extend the range of $\mu$ from $(d/2) + 1 < \mu < d$ to, say, $\sqrt{2d} < \mu < d$, as it was done by Siegel \cite{siegel1} and Dyson \cite{dyson} in the context of (what was later called) the Thue-Siegel principle.

Let us see an application of Theorems \ref{thm:main_theorem_thue} and \ref{thm:main_theorem_thue-mahler}. For an integer $n \geq 3$, let
$$
\Psi_n(x, y) = \prod\limits_{\substack{1 \leq k < \frac{n}{2}\\\gcd(k, n) = 1}}\left(x - 2\cos\left(\frac{2\pi k}{n}\right)y\right)
$$
denote the homogenization of the minimal polynomial of $2\cos\left(\frac{2\pi}{n}\right)$. Then $\Psi_n$ has degree $d = \varphi(n)/2$, where $\varphi(n)$ is the Euler's totient function. Further, the Galois group of $\Psi_n(x, 1)$ has order $d$ \cite[Lemma 3.1]{mosunov1}. Assume that $d \geq 5$ and let $M = |sv - tu|^{-1/2} \cdot \left(\begin{smallmatrix}s & u\\t & v\end{smallmatrix}\right)$ be an element of $\#\Aut'|F|$. Since $d \geq 5$, it follows from Lemma \ref{lem:automorphisms} (see Section \ref{sec:theory}) that there exists a positive integer $j$ coprime to $n$, with $1 \leq j < n/2$, such that
$$
2\cos\left(\frac{2\pi j}{n}\right) = \frac{2\cos\left(\frac{2\pi}{n}\right)v - u}{-2\cos\left(\frac{2\pi}{n}\right) t + s}.
$$
By \cite[Lemma 3.5]{mosunov1}, it must be the case that $s \neq 0$, $s = v$ and $t = u = 0$. Thus, $M = |sv - tu|^{-1/2} \cdot \left(\begin{smallmatrix}s & u\\t & v\end{smallmatrix}\right) = \left(\begin{smallmatrix}\pm 1 & 0\\0 & \pm 1\end{smallmatrix}\right)$. Hence $\Aut'|F| = \left\{\left(\begin{smallmatrix}-1 & 0\\0 & -1\end{smallmatrix}\right), \left(\begin{smallmatrix}1 & 0\\0 & 1\end{smallmatrix}\right)\right\}$, and so the following results hold.

\begin{cor} \label{cor:main_theorem_thue}
Let $n$ be a positive integer such that $\varphi(n) \geq 14$. Let $\lambda$ be such that $0 \leq \lambda < f(d)$, where $f(d)$ is defined in (\ref{eq:bound_on_lambda}). Let $p$ be prime, $k$ a positive integer, and consider the Diophantine equation
\begin{equation} \label{eq:thue_equation_2}
|\Psi_n(x, y)| = tp^k.
\end{equation}
Provided that $p^k$ is sufficiently large, the equation (\ref{eq:thue_equation_2}) has either no solutions in integers $(x, y, t)$ such that
$$
\gcd(x, y) = 1 \quad \text{and} \quad 1 \leq t \leq (p^k)^\lambda,
$$
or exactly two solutions, namely $(x, y, t)$ and $(-x, -y, t)$.
\end{cor}

\begin{cor} \label{cor:main_theorem_thue-mahler}
Let $n$ be a positive integer such that $\varphi(n) \geq 14$. Let $\lambda$ be such that
$$
0 \leq \lambda < 1 - 8.1/(d + 2).
$$
Let $p$ be prime and consider the Diophantine equation
\begin{equation} \label{eq:thue-mahler_equation_2}
|\Psi_n(x, y)| = tp^z.
\end{equation}
Provided that $p$ is sufficiently large, the number of solutions to (\ref{eq:thue-mahler_equation_2}) in integers $(x, y, z, t)$ such that
$$
\gcd(x, y) = 1, \quad z \geq 1 \quad \text{and} \quad 1 \leq t \leq (p^z)^\lambda
$$
is at most
$$
2\left\lfloor 1 + \frac{11.51 + 1.5\log d + \log\left((d - 2.05)/(1 + \lambda)\right)}{\log((d - 2.05)/(1 + \lambda) - 0.5d)}\right\rfloor.
$$
\end{cor}
If we let $d = \varphi(n)/2$ and $\lambda = 1 - 8.5/(d + 2)$, then it is a consequence of Corollary \ref{cor:main_theorem_thue} that the number of solutions in integers $(x, y, z, t)$ to (\ref{eq:thue-mahler_equation_2}) does not exceed $166$ for all $d \geq 7$ and it does not exceed $6$ for all $d \geq 10^{15}$.

The article is structured as follows. In Section \ref{sec:theory} we outline a number of auxiliary results, which are used in later sections. We recommend the reader to skip this section and use it as a reference when reading proofs of Theorems \ref{thm:main_theorem_thue} and \ref{thm:main_theorem_thue-mahler}, which are outlined in Sections \ref{sec:main_theorem_thue} and \ref{sec:main_theorem_thue-mahler}, respectively.

\section{Auxiliary Results} \label{sec:theory}

This section contains several definitions and results, which we utilize in the remaining part of the article. We recommend the reader to skip this section and refer to it when reading the proofs outlined in Sections \ref{sec:main_theorem_thue} and \ref{sec:main_theorem_thue-mahler}.

We begin with a number of definitions. For an arbitrary polynomial $P \in \Z[x_1, x_2, \ldots, x_n]$, we let $H(P)$ denote the maximum of Archimedean absolute values of its coefficients, and refer to this quantity as the \emph{height} of $P$. For a point $(x_1, x_2, \ldots, x_n) \in \Z^n$, we define
$$
H(x_1, x_2, \ldots, x_n) = \max\limits_{i = 1, 2, \ldots, n}\left\{|x_i|\right\}
$$
and refer to this quantity as the \emph{height} of $(x_1, x_2, \ldots, x_n)$.

Let $P \in \mathbb C[x]$ be a polynomial that is not identically equal to zero, with leading coefficient $c_P$. The \emph{Mahler measure} of $P$, denoted $M(P)$, is defined to be $M(P) = |c_P|$ if $P(x)$ is the constant polynomial and
$$
M(P) = |c_P|\prod\limits_{i = 1}^d\max\{1, |\alpha_i|\}
$$
otherwise, where $\alpha_1, \ldots, \alpha_d \in \mathbb C$ are the roots of $P$. For a binary form $Q \in \mathbb C[x, y]$, we define the Mahler measure of $Q$ as $M(Q) = M(Q(x, 1))$. The following lemma is a consequence of a well-known result of Lewis and Mahler \cite{lewis-mahler}. See \cite{mosunov2} for its proof.

\begin{lem}[see {\cite[Lemma 2.6]{mosunov2}}] \label{lem:lewis-mahler}
Let
$$
F(x, y) = c_dx^d + c_{d-1}x^{d-1}y + \cdots + c_0y^d
$$
be a binary form of degree $d \geq 2$ with integer coefficients such that $c_0c_d \neq 0$. Let $x_0, y_0$ be nonzero integers. There exists a root $\alpha$ of $F(x, 1)$ such that
$$
\min\left\{\left|\alpha - \frac{x_0}{y_0}\right|, \left|\alpha^{-1} - \frac{y_0}{x_0}\right|\right\} \leq \frac{C|F(x_0, y_0)|}{H(x_0, y_0)^d},
$$
where
$$
C =  \frac{2^{d-1}d^{(d-1)/2}M(F)^{d-2}}{|D(F)|^{1/2}}.
$$
\end{lem}

\begin{lem}[Thunder, {\cite[Lemma 2]{thunder2}}] \label{lem:thunders_lemma}
Let $p$ be a rational prime and let $\overline{\Q_p}$ denote the algebraic closure of the field of $p$-aidc numbers $\Q_p$. Let $F \in \Z[x, y]$ be an irreducible homogeneous polynomial of degree $d \geq 2$ and content one, and denote the roots of $F(x, 1)$ by $\alpha_1, \ldots, \alpha_d \in \overline{\Q_p}$. Let $x$ and $y$ be coprime integers. If $i_0$ is an index with
$$
\frac{|x - \alpha_{i_0} y|_p}{\max\{1, |\alpha_{i_0}|_p\}} = \min\limits_{1 \leq i \leq d}\left\{\frac{|x - \alpha_i y|_p}{\max\{1, |\alpha_i|_p\}}\right\},
$$
then
$$
\frac{|x - \alpha_{i_0}y|_p}{\max\{1, |\alpha_{i_0}|_p\}} \leq \frac{|F(x, y)|_p}{|D(F)|_p^{1/2}}.
$$
Further, if $|F(x, y)|_p < |D(F)|_p^{1/2}$, then the index $i_0$ above is unique and $\alpha_{i_0} \in \Q_p$.
\end{lem}

The following three results were established in \cite{mosunov2}. Lemma \ref{lem:non-archimedean_gap_principle} states the generalized non-Archimedean gap principle, which plays a crucial role in the proof of Theorem \ref{thm:main_theorem_thue}. In turn, Lemma \ref{lem:count_large_approximations} follows from the combination of the generalized gap principle and the Thue-Siegel principle, as formulated by Bombieri and Mueller in \cite[Section II]{bombieri-mueller}.

\begin{lem}[see {\cite[Theorem 1.2]{mosunov2}}] \label{lem:non-archimedean_gap_principle} Let $p$ be a rational prime. Let $\alpha \in \Q_p$ be a $p$-adic algebraic number of degree $d \geq 3$ over $\Q$ and let $\beta$ be irrational and in $\Q(\alpha)$. Let $\mu$ be a real number such that $(d/2) + 1 < \mu < d$ and let $C_0$ be a positive real number. There exist positive real numbers $C_{1}$ and $C_{2}$, that are explicitly computable in terms of $\alpha$, $\beta$, $\mu$ and $C_0$, with the following property. If $x_1/y_1$ and $x_2/y_2$ are rational numbers in lowest terms such that $H(x_2, y_2) \geq H(x_1, y_1) \geq C_{1}$ and
$$
\left|y_1\alpha - x_1\right|_p < \frac{C_0}{H(x_1, y_1)^\mu}, \quad \left|y_2\beta - x_2\right|_p < \frac{C_0}{H(x_2, y_2)^\mu},
$$
then at least one of the following holds:
\begin{itemize}
\item  $H(x_2, y_2) > C_{2}^{-1} H(x_1, y_1)^{\mu - d/2}$;
 \item There exist integers $s, t, u, v$, with $sv - tu \neq 0$, such that
$$
\beta = \frac{s\alpha + t}{u\alpha + v} \quad \text{and} \quad \frac{x_2}{y_2} = \frac{sx_1 + ty_1}{ux_1 + vy_1}.
$$
\end{itemize}
\end{lem}

\begin{lem}[see {\cite[Proposition 7.3]{mosunov2}}] \label{lem:automorphisms}
Let $F \in \Z[x, y]$ be an irreducible binary form of degree $d \geq 3$ and let $c_d$ denote the coefficient of $x^d$ in $F$. Let $\alpha_1, \ldots, \alpha_d$ be the roots of $F(x, 1)$. There exists an index $j \in \{1, \ldots, d\}$ such that
$$
\alpha_j = \frac{v\alpha_1 - u}{-t\alpha_1 + s}
$$
for some integers $s$, $t$, $u$ and $v$ if and only if the matrix
$$
M = \frac{1}{\sqrt{|sv - tu|}}
\begin{pmatrix}
s & u\\
t & v
\end{pmatrix}
$$
is in $\Aut'|F|$. Furthermore, if $M \in \Aut'|F|$, then $|sv - tu| = \left|\frac{F(s, t)}{c_d}\right|^{2/d}$.
\end{lem}

For an irrational number $\alpha$, the \emph{orbit} of $\alpha$ is the set
$$
\orb(\alpha) = \left\{\frac{v\alpha - u}{-t\alpha + s}\ \colon\ s,t,u,v \in \Z,\ sv-tu \neq 0\right\}.
$$

\begin{lem} \label{lem:count_large_approximations}
Let $K = \R$ or $\Q_p$, where $p$ is a rational prime, and denote the standard absolute value on $K$ by $|\quad|$. Let $\alpha_1 \in K$ be an algebraic number of degree $d \geq 3$ over $\Q$ and $\alpha_2, \alpha_3, \ldots, \alpha_n$ be distinct elements of $\Q(\alpha_1)$, different from $\alpha_1$, each of \mbox{degree $d$.} Let $\mu$ be such that $(d/2) + 1 < \mu < d$. Let $C_0$ be a real number such that $C_0 > (4e^A)^{-1}$, where
\begin{equation} \label{eq:A}
A = 500^2\left(\log \max\limits_{i = 1, \ldots, n}\{M(\alpha_i)\} + \frac{d}{2}\right).
\end{equation}
There exists a positive real number $C_3$, which is explicitly computable in terms of $\alpha_1, \alpha_2, \ldots, \alpha_n$, $\mu$ and $C_0$, with the following property. The total number of rationals $x/y$ in lowest terms, which satisfy \mbox{$H(x, y) \geq C_3$} and
\begin{equation} \label{eq:roths_inequality}
\left|\alpha_j - \frac{x}{y}\right| < \frac{C_0}{H(x, y)^\mu}
\end{equation}
for some $j \in \{1, 2, \ldots, n\}$ is less than
$$
\gamma\left\lfloor1 + \frac{11.51 + 1.5\log d + \log \mu}{\log(\mu - d/2)}\right\rfloor,
$$
where
\begin{equation} \label{eq:gamma}
\gamma = \max\{\gamma_1, \ldots, \gamma_n\}, \quad \gamma_i = \#\{j \colon \alpha_j \in \orb(\alpha_i)\}.
\end{equation}
\end{lem}

Notice that when degree $d$ extension $\mathbb Q(\alpha)/\mathbb Q$ is Galois and $\alpha = \alpha_1, \ldots, \alpha_d$ are the algebraic conjugates of $\alpha$, then for $d \geq 3$ it follows from Lemma \ref{lem:automorphisms} that $\alpha_j = (v\alpha - u)/(-t\alpha + s)$ if and only if $M = |sv - tu|^{-1/2}\cdot \left(\begin{smallmatrix}s & u\\t & v\end{smallmatrix}\right)$ is an element of $\Aut'|F|$. Thus, in this case, the quantity $\gamma$ in (\ref{eq:gamma}) does not exceed $\#\Aut'|F|$. This fact plays an important role in the proof of Theorem \ref{thm:main_theorem_thue-mahler}.

\section{Proof of Theorem \ref{thm:main_theorem_thue}} \label{sec:main_theorem_thue}
By Roth's Theorem \cite{roth}, for every complex root $\alpha$ of $F(x, 1)$ there exist only finitely many nonzero integers $x, y$ such that $\min\left\{|\alpha - x/y|, \left|\alpha^{-1} - y/x\right|\right\} \leq H(x, y)^{-2.05}$. Since
$$
|F(x, y)| \leq (d + 1)H(F)H(x, y)^d
$$
and $|F(x, y)| = tp^k$, we have
$$
\frac{tp^k}{(d + 1)H(F)} \leq H(x, y)^d.
$$
Hence by choosing a large enough $p^k$ we can increase $H(x, y)$ and make it so large that the inequality $\min\left\{|\alpha - x/y|, \left|\alpha^{-1} - y/x\right|\right\} \leq H(x, y)^{-2.05}$ is no longer satisfied for every complex root $\alpha$ of $F(x, 1)$.

Define
$$
C_0 = \frac{2^{d-1}d^{(d-1)/2}M(F)^{d-2}}{|D(F)|^{1/2}}.
$$
Assume that there exists a solution $(x, y, t)$ to (\ref{eq:thue_equation_1}). By Lemma \ref{lem:lewis-mahler},
$$
\min\left\{\left|\alpha - \frac{x}{y}\right|, \left|\alpha^{-1} - \frac{y}{x}\right|\right\} \leq \frac{C_0tp^k}{H(x, y)^d}.
$$
From our choice of $p^k$ and the above inequality it follows that
$$
\frac{1}{H(x, y)^{2.05}} < \min\left\{\left|\alpha - \frac{x}{y}\right|, \left|\alpha^{-1} - \frac{y}{x}\right|\right\} \leq \frac{C_0tp^k}{H(x, y)^d},
$$
which is equivalent to
\begin{equation} \label{eq:bound_on_H}
H(x, y) < (C_0tp^k)^{1/(d - 2.05)}.
\end{equation}

Since $t \leq (p^k)^\lambda$,
$$
tp^k \leq (p^k)^{1 + \lambda} \leq |F(x, y)|_p^{-(1 + \lambda)}.
$$
Combining this inequality with (\ref{eq:bound_on_H}), we get
$$
H(x, y)^{d - 2.05}
< C_0tp^k \leq C_0|F(x, y)|_p^{-(1 + \lambda)}.\\
$$
We conclude that
\begin{equation} \label{eq:to_be_applied}
|F(x, y)|_p < \frac{C_0^{1/(1 + \lambda)}}{H(x, y)^\mu},
\end{equation}
where
$$
\mu = \frac{d - 2.05}{1 + \lambda}.
$$

Next, we take $p^k$ sufficiently large that
$$
p^k > |D(F)|.
$$
Then
$$
|F(x, y)|_p \leq p^{-k} < |D(F)|^{-1} \leq |D(F)|_p.
$$
By Lemma \ref{lem:thunders_lemma} there exists a unique $p$-adic root $\alpha \in \Q_p$ of $F(x, 1)$ such that
$$
\frac{|y\alpha - x|_p}{\max\{1, |\alpha|_p\}} \leq \frac{|F(x, y)|_p}{|D(F)|_p^{1/2}}.
$$
Let $c_d$ denote the coefficient of $x^d$ in $F$. Since $c_d\alpha$ is an algebraic integer, we see that $|c_d\alpha|_p \leq 1$, so $\max\{1, |\alpha|_p\} \leq |c_d|_p^{-1}$. Combining this inequality with (\ref{eq:to_be_applied}), we obtain
\begin{align*}
|y\alpha - x|_p
& < \frac{\max\{1, |\alpha|_p\}}{|D(F)|_p^{1/2}}|F(x, y)|_p\\
& \leq \frac{C_1}{H(x, y)^\mu},
\end{align*}
where
$$
C_1 = C_0^{1/(1 + \lambda)}c_d|D(F)|^{1/2}.
$$

Now, assume that there exist two solutions $(x_1, y_1, t_1)$ amd $(x_2, y_2, t_2)$ to (\ref{eq:thue_equation_1}), ordered so that $H(x_2, y_2) \geq H(x_1, y_1)$. Then it follows from the discussion above that there exist $p$-adic roots $\alpha, \beta \in \Q_p$ of $F(x, 1)$ such that
$$
|y_1\alpha - x_1|_p < \frac{C_1}{H(x_1, y_1)^\mu}, \quad |y_2\beta - x_2|_p < \frac{C_1}{H(x_2, y_2)^\mu}.
$$
Since $(d/2) + 1 < \mu < d$, it follows from Lemma \ref{lem:non-archimedean_gap_principle} that there exists positive numbers $C_2$ and $C_3$, which depend on $C_1$, $\mu$, and $F$, but not on $p$, such that if $H(x_2, y_2) \geq H(x_1, y_1) \geq C_2$, then either $H(x_2, y_2) > C_3H(x_1, y_1)^{\mu - d/2}$, or $\alpha, \beta$ and $x_1/y_1, x_2/y_2$ are connected by means of a linear fractional transformation, or both. By choosing $p^k$ sufficiently large we can always ensure that $H(x_1, y_1) \geq C_2$. We obtain an upper bound on $H(x_2, y_2)$ by combining (\ref{eq:bound_on_H}) with the inequality $|F(x_1, y_1)| \leq (d + 1)H(F)H(x_1, y_1)^d$:
\begin{align*}
H(x_2, y_2)
& < (C_0t_2p^k)^{1/(d - 2.05)}\\
& \leq \left(C_0(p^k)^{1 + \lambda}\right)^{1/(d - 2.05)}\\
& \leq \left(C_0(t_1p^k)^{1 + \lambda}\right)^{1/(d - 2.05)}\\
& = \left(C_0|F(x_1, y_1)|^{1 + \lambda}\right)^{1/(d - 2.05)}\\
& \leq \left(C_0\left((d + 1)H(F)H(x_1, y_1)^d\right)^{1 + \lambda}\right)^{1/(d - 2.05)}. 
\end{align*}
Merging the above upper bound with the lower bound $H(x_2, y_2) > C_3H(x_1, y_1)^{\mu - d/2}$ results in the inequality
$$
C_3H(x_1, y_1)^{\mu - d/2 - d(1 + \lambda)/(d - 2.05)} < \left(C_0\left((d + 1)H(F)\right)^{1 + \lambda}\right)^{1/(d - 2.05)}.
$$
From our choice of $\lambda$ it follows that the exponent of $H(x_1, y_1)$ is positive, and so $H(x_1, y_1)$ is bounded. Thus, by making $p^k$ (and therefore $H(x_1, y_1)$) sufficiently large we can always ensure that the inequality $H(x_2, y_2) > C_3H(x_1, y_1)^{\mu - d/2}$ does not hold. Then $\alpha, \beta$ and $x_1/y_1, x_2/y_2$ are connected by means of a linear fractional transformation:
$$
\beta = \frac{v\alpha - u}{-t\alpha + s} \quad \text{and} \quad \frac{x_2}{y_2} = \frac{vx_1 - uy_1}{-tx_1 + sy_1},
$$
where $s, t, u, v \in \Z$ and $sv - tu \neq 0$. By Lemma \ref{lem:automorphisms}, the matrix
$$
M = \frac{1}{\sqrt{|sv - tu|}}
\begin{pmatrix}
s & u\\
t & v
\end{pmatrix}
$$
is an element of $\Aut' |F|$. Hence the number of solutions $(x, y, t)$ to (\ref{eq:thue_equation_1}) is at most $\#\Aut' |F|$.

\section{Proof of Theorem \ref{thm:main_theorem_thue-mahler}} \label{sec:main_theorem_thue-mahler}
The beginning of the proof is similar to the proof of Theorem \ref{thm:main_theorem_thue}. By Roth's Theorem \cite{roth}, for every root $\alpha$ of $F(x, 1)$ there exist only finitely many nonzero integers $x, y$ such that $\min\left\{|\alpha - x/y|, |\alpha^{-1} - y/x|\right\} \leq H(x, y)^{-2.05}$. Since $|F(x, y)| \leq (d + 1)H(F)H(x, y)^d$ and $|F(x, y)| = tp^z$, we have
$$
\frac{p}{(d + 1)H(F)} \leq \frac{tp^z}{(d + 1)H(F)} = \frac{|F(x, y)|}{(d + 1)H(F)} \leq H(x, y)^d.
$$
Hence by choosing a large enough $p$ we can increase $H(x, y)$ and make it so large that the inequality $\min\left\{|\alpha - x/y|, |\alpha^{-1} - y/x|\right\} \leq H(x, y)^{-2.05}$ is no longer satisfied for every complex root $\alpha$ of $F(x, 1)$.

Now, assume that there exists a solution $(x, y, z, t)$ of (\ref{eq:thue-mahler_equation_1}). As in the proof of Theorem \ref{thm:main_theorem_thue}, for our choice of $p$ the inequality
\begin{equation} \label{eq:bound_on_H_2}
H(x, y) < (C_0tp^z)^{1/(d - 2.05)}
\end{equation}
holds, where
$$
C_0 = \frac{2^{d-1}d^{(d-1)/2}M(F)^{d-2}}{|D(F)|^{1/2}}.
$$
Since
$$
tp^z \leq (p^z)^{1 + \lambda} \leq |F(x, y)|_p^{-(1 + \lambda)},
$$
it follows from (\ref{eq:bound_on_H_2}) that
$$
|F(x, y)|_p < \frac{C_0^{1/(1 + \lambda)}}{H(x, y)^{\mu}},
$$
where
$$
\mu = \frac{d - 2.05}{1 + \lambda}.
$$

We take $p$ sufficiently large that
$$
p > |D(F)|.
$$
Then
$$
|F(x, y)|_p \leq p^{-1} < |D(F)|^{-1} \leq |D(F)|_p.
$$
Let $c_d$ denote the coefficient of $x^d$ in $F$. By Lemma \ref{lem:thunders_lemma} there exists a unique $p$-adic root $\alpha \in \Q_p$ of $F(x, 1)$ such that
$$
|y\alpha - x|_p \leq \frac{\max\{1, |\alpha|_p\}}{|D(F)|_p^{1/2}}|F(x, y)|_p < \frac{C_1}{H(x, y)^\mu},
$$
where
$$
C_1 = C_0^{1/(1 + \lambda)}c_d|D(F)|^{1/2}.
$$
Note that $C_1$ is independent \mbox{of $p$.} Further, we can ensure that $p \nmid y$ by adjusting our choice of $p$ as follows:
$$
p > c_d.
$$
Indeed, if $p \mid y$, then $p$ does not divide $x$, because $x$ and $y$ are coprime. Since $z \geq 1$, it is evident from equation
$$
c_dx^d + y(c_{d-1}x^{d-1} + \cdots + c_0y^{d-1}) = \pm tp^z
$$
that $p$ divides $c_d$, in contradiction to our choice of $p$. Then $|y|_p = 1$, and so for any $\alpha \in \Q_p$ we have
$$
\left|\alpha - \frac{x}{y}\right|_p = |y\alpha - x|_p.
$$
Therefore
$$
\left|\alpha - \frac{x}{y}\right|_p < \frac{C_1}{H(x, y)^\mu}.
$$

Let $\alpha_1, \alpha_2, \ldots, \alpha_d$ be the roots of $F(x, 1)$. Since $(d/2) + 1 < \mu < d$, we apply Lemma \ref{lem:count_large_approximations} and conclude that there exists a positive number $C_2$, which depends on $C_1, \mu, \alpha_1, \alpha_2, \ldots, \alpha_d$, but not on $p$, such that the number of rationals $x/y$ in lowest terms satisfying $H(x, y) \geq C_2$ and
\begin{equation} \label{eq:roths_inequality_2}
\left|\alpha_j - \frac{x}{y}\right|_p < \frac{C_1}{H(x, y)^\mu}
\end{equation}
for some $j \in \{1, 2, \ldots, d\}$ is less than
$$
\#\Aut' |F|\cdot \left\lfloor 1 + \frac{11.51 + 1.5\log d + \log \mu}{\log(\mu - 0.5d)}\right\rfloor.
$$
If we choose $p$ so that $p \geq (d + 1)H(F)C_2^d$, then
$$
C_2^d \leq \frac{p}{(d + 1)H(F)} \leq \frac{tp^z}{(d + 1)H(F)} = \frac{|F(x, y)|}{(d + 1)H(F)} \leq H(x, y)^d,
$$
so the inequality $H(x, y) \geq C_2$ is satisfied. Since all solutions $(x, y, z, t)$ to (\ref{eq:thue-mahler_equation_1}), including those that satisfy $H(x, y) \geq C_2$, also satisfy (\ref{eq:roths_inequality_2}), the result follows.

\section*{Acknowledgements}

The author is grateful to Prof.\ Cameron L.\ Stewart for his wise supervision.

\bibliography{mosunov-number-of-solutions}

\begin{thebibliography}{10}

\bibitem{bombieri-mueller}
E.~Bombieri and J.~Mueller.
\newblock On effective measures of irrationality for $\sqrt[r]{\frac{a}{b}}$
  and related numbers.
\newblock {\em J. Reine Angew. Math.}, 342:173--196, 1983.

\bibitem{bombieri-schmidt}
E.~Bombieri and W.~M. Schmidt.
\newblock On {T}hue's equation.
\newblock {\em Invent. Math.}, 88:69--81, 1987.

\bibitem{dyson}
F.~Dyson.
\newblock The approximation of algebraic numbers by rationals.
\newblock {\em Acta Math.}, 79:225--240, 1947.

\bibitem{erdos-mahler}
P.~Erd\H{o}s and K.~Mahler.
\newblock On the number of integers which can be represented by a binary form.
\newblock {\em J. London Math. Soc.}, 13:134--139, 1938.

\bibitem{erdos-stewart-tijdeman}
P.~Erd\H{o}s, C.~L. Stewart, and R.~Tijdeman.
\newblock Some {D}iophantine equations with many solutions.
\newblock {\em Comp. Math.}, 66:37--56, 1988.

\bibitem{evertse2}
J.-H. Evertse.
\newblock On equations in {$S$}-units and the {T}hue-{M}ahler equation.
\newblock {\em Invent. Math.}, 75:561--584, 1984.

\bibitem{lewis-mahler}
D.~Lewis and K.~Mahler.
\newblock Representation of integers by binary forms.
\newblock {\em Acta Arith.}, 6:333--363, 1961.

\bibitem{mahler2}
K.~Mahler.
\newblock Zur {A}pproximation algebraischer {Z}ahlen. {II}. \"{U}ber die
  {A}nzahl der {D}arstellungen ganzer {Z}ahlen durch {B}in\"arformen.
\newblock {\em Math. Ann.}, 108:37--55, 1933.

\bibitem{mosunov1}
A.~Mosunov.
\newblock On the automorphism group of a binary form associated with algebraic
  grigonometric quantities.
\newblock {\em J. of Number Theory}, 2022.

\bibitem{mosunov2}
A.~Mosunov.
\newblock On the generalization of the gap principle.
\newblock {\em Period. Math. Hungar. (to appear)}, 2022.

\bibitem{roth}
K.~Roth.
\newblock Rational approximations to algebraic numbers.
\newblock {\em Mathematika}, 2(3):1--20, 1955.

\bibitem{siegel1}
C.~L. Siegel.
\newblock {A}pproximation algebraischer {Z}ahlen.
\newblock {\em Math. Zeit.}, 10:173--213, 1921.

\bibitem{stewart2}
C.~L. Stewart.
\newblock On the number of solutions of polynomial congruences and {T}hue
  equations.
\newblock {\em J. Amer. Math. Soc.}, 4:793--835, 1991.

\bibitem{stewart3}
C.~L. Stewart.
\newblock Cubic {T}hue equations with many solutions.
\newblock {\em International Mathematics Research Notices}, page rnn040, 2008.

\bibitem{thue}
A.~Thue.
\newblock \"{U}ber {A}nn\"aherungswerte algebraischer {Z}ahlen.
\newblock {\em J. Reine Angew. Math.}, 135:284--305, 1909.

\bibitem{thunder2}
J.~L. Thunder.
\newblock Thue equations and lattices.
\newblock {\em Illinois J. Math.}, 59(4):999--1023, 2015.

\end{thebibliography}
\bibliographystyle{plain}
\end{document}